\documentclass[12pt]{article}
\usepackage{mathrsfs}
\usepackage{dsfont}
\usepackage{amsthm}
\usepackage{mathrsfs}
\usepackage{amsmath}
\usepackage{amsfonts}
\usepackage[colorlinks,linkcolor=black,anchorcolor=blue,citecolor=blue]{hyperref}
\usepackage{amssymb, amsmath, cite}

\newtheorem{theorem}{Theorem}[section]

\newcommand\be{\begin{equation}}
\newcommand\ee{\end{equation}}
\newcommand\ber{\begin{eqnarray}}
\newcommand\eer{\end{eqnarray}}
\newcommand\berr{\begin{eqnarray*}}
\newcommand\eerr{\end{eqnarray*}}
\newcommand\bea{\begin{eqnarray}}
\newcommand\eea{\end{eqnarray}}

\newcommand{\nn}{\nonumber}

\newcommand{\dd}{\mathrm{d}}
\newcommand\e{\mathrm{e}}\newcommand\pa{\partial}
{\setlength\arraycolsep{2pt}

\begin{document}

\title{Existence of non-Abelian vortices in a coupled 4D--2D quantum field theory}
\author{Yilu Xu\footnote{E-mail address:
104753180630@vip.henu.edu.cn.}, Shouxin Chen\\School of Mathematics and Statistics, Henan University\\
Kaifeng, Henan 475004, PR China}
\date{}
\maketitle

\begin{abstract}
Vortices produce locally concentrated field configurations and are solutions to the nonlinear partial differential equations systems of complicated structures. In this paper, we establish the existence and uniqueness for solutions of the gauged non--Abelian vortices in a coupled 4D--2D quantum field theory by researching the nonlinear elliptic equations systems with exponential terms in $\mathbb{R}^{2}$ using the calculus of variations. In addition, we obtain the asymptotic behavior of the solutions at infinity and the quantized integrals in $\mathbb{R}^{2}$.
\end{abstract}

\begin{enumerate}

\item[]

\item[]
{MSC numbers(2020):} 35B38, 35B40, 81V25

\end{enumerate}

\section{Introduction}\label{s1}
\setcounter{equation}{0}

Gauge fields originated from classical electromagnetism and have become the cornerstone of modern physics after Yang and Mills\cite{YM} generalized to non--Abelian. The study of the partial differential equation problems arising in field theory is of great importance for research areas at the intersection of fundamental mathematics and mathematical physics. Depending on the spatial dimensions, the static topological solutions for the gauge field equations can be categorized into four types of solutions, namely, instantons, monopoles, vortices and domain walls, which are collectively referred to as solitons\cite{Manton}. Solitons play a fundamental role in gauge field theory, including quantum field theory and classical field theory and are crucial for describing a variety of fundamental interactions and rich phenomena\cite{Manton}.

Vortices, as two--dimensional static solutions of the gauge field equations, have always attracted attention and interest due to their essential applications in various fundamental areas of theoretical physics. While describing physical theories in different contexts, vortices usually play very significant roles in understanding key phenomena. In the Ginzburg--Landau superconductivity theory, magnetically charged vortices, the simplest Abelian vortices, were first discovered in 1957 by Abrikosov\cite{Abrokosov} in the form of a typical mixed state in a type II superconductor\cite{Jaffe}. In 1973, Nielsen and Olesen\cite{Nielsen} demonstrated that vortices exist in the Abelian--Higgs model in the context of quantum field theory. Furthermore, according to the Julia--Zee theorem\cite{Julia}, the vortex equations in the classical Abelian Higgs model are known as the Ginzburg--Landau equations in the static limit. In high--energy physics, cosmic strings\cite{Hindmarsh} are the gravitational vortices created when the Einstein gravity theory is coupled to the Abelian--Higgs model, which contributed to the formation of matter in the early universe\cite{Abel,Vilenkin}. Subsequently, non--Abelian vortices which generally arise in the color--flavor locked\cite{Shifman} phases of the non--Abelian gauge theory have attracted lots of attention as a consequence of their meaningful roles in grand unified theories. In the unified theory of the Glasgow--Weinberg--Salam weak forces and the electromagnetic forces, vortex solutions are generated by W and Z particle fields\cite{Ambjorn} and lead to the anti--Meissner effect\cite{Bartolucci,YangBook}. Currently, the electroweak theory of Glashow--Weinberg--Salam\cite{Actor} is not only a vital unified field theory available to describe fundamental interactions, but also one of the most successful non--Abelian gauge field theories. One of the fundamental problems in particle physics is that of quark confinement\cite{Mandelstam,Shifman} which is the inability of the quarks that make up elementary particles to be observed in isolation. The vortices generated by non--Abelian gauge fields can make the attraction between quarks and anti--quarks independent of their distance, thus achieving linear confinement between quarks, theoretical physicists have conducted extensive research on them and derived a wide range of nonlinear vortex equations with rich features\cite{Hanany,Tong}. In addition, both electrically and magnetically charged vortices have applications in a vast range of areas in condensed matter physics such as the quantum Hall effect\cite{Sokoloff}, high temperature superconductors\cite{Khomskii}, the Bose--Einstein condensates\cite{Kawaguchi}, holographic superconductors\cite{Hartnoll,Tallarita}, superfluids\cite{Shevchenko}, optics\cite{Bezryadina} and quantum chromodynamics\cite{Fujimoto}.

It is worth noting that we will refer to vortices generated by systems where the gauge fields are coupled to part or all of the color--flavor diagonal global symmetry as gauged non--Abelian vortices in this paper. The purpose of the present paper is to establish the existence and uniqueness theorem for the gauged non--Abelian vortex model in a coupled 4D--2D quantum field theory proposed in the work of Bolognesi et al.\cite{Bo}. Inspired by the methods of \cite{Chen0,Chen1,Chen2,Lieb,YangBook}, we obtain existence and uniqueness results for non--Abelian vortex solutions and build the explicit decay estimates for planar solutions.

We complete this section with a brief overview of the paper, explaining how the above results are organized. The paper is setup as follows. In Section 2, we introduce a system of nonlinear ordinary differential equations from the gauged non--Abelian vortex model proposed in \cite{Bo}. In Section 3, we derive the governing system which is a nonlinear elliptic equations system with exponential terms from the model in Section 2, and present the existence and uniqueness theorem over the full plane. In Section 4, we establish the existence and uniqueness result for solutions of the vortex equations over the full plane, applying the calculus of variations first developed by Jaffe and Taubes\cite{Jaffe} for the Abelian Higgs model and the convexity of an action functional. Furthermore, we demonstrate the exponential decay properties of the solutions and their derivative at infinity and utilize them to acquire the anticipated quantized integrals in $\mathbb{R}^{2}$.

\section{Gauged non--Abelian vortex model}\label{s2}
\setcounter{equation}{0}

In this section, we start by a review of the non--Abelian vortex model with gauge group derived by Bolognesi et al.\cite{Bo}. Here we will give a rough description, as details can be found in \cite{Bo}. In the classical Abelian Higgs model, one starts from a complex scalar field $Q$ that lies in the fundamental or defining representation of $U(1)$. Bolognesi et al. study the extension of this theory into the situation that the gauge group of the type
\be\label{2.1}
G=U_{0}(1)\times G_{L}\times G_{R},
\ee
where $G_{L}=G_{R}=SU(N)$. The matter sector is composed of a complex scalar field $Q$ in the dual fundamental representation of the two $SU(N)$ factors with unit charge with respect to $U_{0}(1)$. Following Bolognesi et al.\cite{Bo}, we consider the truncated bosonic sector of the $\mathcal{N}=2$ supersymmetric theory BPS--saturated action as
\ber\label{2.2}
&&\mathcal{L}=-\frac{1}{2}\mathrm{Tr}(F_{\mu\nu}^{(l)}F^{(l)\mu\nu})-\frac{1}{2}\mathrm{Tr}(F_{\mu\nu}^{(r)}F^{(r)\mu\nu})-\frac{1}{4}f_{\mu\nu}f^{\mu\nu}+\mathrm{Tr}(D_{\mu}Q^{\dag}D^{\mu}Q)\nn\\
&&~~~~\,~-\frac{g_{0}^{2}}{2}(\mathrm{Tr}Q^{\dag}Q-v_{0}^{2})^{2}-\frac{g_{l}^{2}}{2}(\mathrm{Tr}t^{a}QQ^{\dag})^{2}-\frac{g_{r}^{2}}{2}(\mathrm{Tr}t^{a}Q^{\dag}Q)^{2},
\eer
where $v_{0}^{2}=N\xi$ and $\dag$ denotes the Hermitian conjugate. The covariant derivative is defined by
\be\label{2.3}
D_{\mu}Q=\pa_{\mu}Q-\mathrm{i}g_{l}A_{\mu}^{(l)}Q-\mathrm{i}g_{0}a_{\mu}Q+\mathrm{i}g_{r}QA_{\mu}^{(r)}.
\ee
In the vacuum, the scalar--field condensate takes the form
\be\label{2.4}
\langle Q\rangle=\sqrt{\xi}\mathbf{1}_{N},
\ee
leaving a diagonal $SU(N)$ gauge group unbroken. The fields
\be\label{2.5}
\mathcal{A}_{\mu}=\frac{1}{\sqrt{g_{r}^{2}+g_{l}^{2}}}\left(g_{r}A_{\mu}^{(l)}+g_{l}A_{\mu}^{(r)}\right)
\ee
keep massless in the bulk, however the orthogonal combination
\be\label{2.6}
\mathcal{B}_{\mu}=\frac{1}{\sqrt{g_{r}^{2}+g_{l}^{2}}}\left(g_{r}A_{\mu}^{(l)}-g_{l}A_{\mu}^{(r)}\right)
\ee
and the $U(1)$ field $a_{\mu}$ become massive. The nontrivial first homotopy group
\be\label{2.7}
\pi_{1}\left(\frac{U_{0}(1)\times SU_{L}(N)\times SU_{R}(N)}{SU_{L+R}(N)}\right)=\mathbb{Z},
\ee
implies that the system allows stable vortices.

We are interested in stable vortices. Now, we assume that all the field configurations depend on the transverse coordinates $x$ and $y$ only. Therefore, the vortex solutions can be derived by the BPS completion of the expression for the tension:
\ber\label{2.8}
&&T=\int\mathrm{d}^{2}x\Bigg\{\frac{1}{2}\left(f_{12}+g_{0}\left(\mathrm{Tr}Q^{\dag} Q-N\xi\right)\right)^{2}+\mathrm{Tr}\bigg[\left(F_{12}^{(r)}-g_{r}t^{a}\mathrm{Tr}\left(Q^{\dag} Qt^{a}\right)\right)^{2}\nn\\
&&~~~~~~~~+\left(F_{12}^{(l)}+g_{l}t^{a}\mathrm{Tr}\left(Q^{\dag} t^{a}Q\right)\right)^{2}\bigg]+\mathrm{Tr}|D_{1}Q+\mathrm{i}D_{2}Q|^{2}+g_{0}N\xi f_{12}\Bigg\}.
\eer
According to \eqref{2.8}, we arrive at the BPS equations
\ber
D_{1}Q+\mathrm{i}D_{2}Q&=&0,\label{2.9}\\
f_{12}+g_{0}\left(\mathrm{Tr}Q\dag Q-N\xi\right)&=&0,\label{2.10}\\
F_{12}^{(r)}-g_{r}t^{a}\mathrm{Tr}\left(Q\dag Qt^{a}\right)&=&0,\label{2.11}\\
F_{12}^{(l)}+g_{l}t^{a}\mathrm{Tr}\left(Q\dag t^{a}Q\right)&=&0.\label{2.12}
\eer
For a minimal vortex with a fixed orientation in color--flavor, for example $(1,\mathbf{1}_{N-1})$, we can choose the radially symmetric ansatz in the scalar field
\ber\label{2.13}
Q=\left(
\begin{matrix}
\mathrm{e}^{\mathrm{i}\theta}Q_{1}(r)&0\\[8pt]
0&Q_{2}(r)\mathbf{1}_{N-1}
\end{matrix}
\right),
\eer
while the Abelian and non--Abelian gauge fields can be expressed as the diagonal form
\ber
a_{i}&=&-\frac{1}{g_{0}}\frac{\epsilon_{ij}x_{j}}{r^{2}}\frac{1-f}{N},\label{2.14}\\
A_{i}^{(l)}&=&-\frac{g_{l}}{g'^{2}}\frac{\epsilon_{ij}x_{j}}{r^{2}}\frac{1-f_{NA}}{NC_{N}}T_{N^{2}-1},\label{2.15}\\
A_{i}^{(r)}&=&\frac{g_{r}}{g'^{2}}\frac{\epsilon_{ij}x_{j}}{r^{2}}\frac{1-f_{NA}}{NC_{N}}T_{N^{2}-1},\label{2.16}
\eer
where
\ber\label{2.17}
T_{N^{2}-1}\equiv C_{N}\left(
\begin{matrix}
N-1&0\\[8pt]
0&-\mathbf{1}_{N-1}
\end{matrix}
\right),~~~C_{N}\equiv\frac{1}{\sqrt{2N(N-1)}},~~~g'\equiv\sqrt{g_{l}^{2}+g_{r}^{2}}.
\eer
From the BPS equation \eqref{2.9}--\eqref{2.12} and the above ansatz, setting $r=|x|$, then we can get that the profile functions satisfy
\ber
\frac{f'}{r}-g_{0}^{2}N\left[Q_{1}^{2}+(N-1)Q_{2}^{2}-N\xi\right]&=&0,\label{2.18}\\
\frac{f'_{NA}}{r}-g'^{2}\frac{Q_{1}^{2}-Q_{2}^{2}}{2}&=&0,\label{2.19}\\
rQ'_{1}-Q_{1}\left(\frac{(N-1)f_{NA}+f}{N}\right)&=&0,\label{2.20}\\
rQ'_{2}-Q_{2}\left(\frac{-f_{NA}+f}{N}\right)&=&0\label{2.21}
\eer
and the boundary conditions
\ber
f(r)&=&1,~~~f_{NA}(r)=1,~~~\text{as}~~r\to 0,\label{2.22}\\
f(r)&=&0,~~~f_{NA}(r)=0,~~~Q_{1}(r)=\sqrt{\xi},~~~Q_{2}(r)=\sqrt{\xi},~~~\text{as}~~r\to \infty.\label{2.23}
\eer
It should be noted that the finite vortex points of both $Q_{1}$ and $Q_{2}$ are concentrated at the origin in the radial symmetry situation, that is to say, $Q_{1}(r)=Q_{2}(r)=0$ as $r\to 0$. Moreover, in view of the equations \eqref{2.18}--\eqref{2.21}, we see that $f$ is not equal to $f_{NA}$ and $Q_{1}$ is not equal to $Q_{2}$. For the nonlinear ordinary differential equations \eqref{2.18}--\eqref{2.21} subject to the boundary conditions \eqref{2.22}--\eqref{2.23}, we know that the existence of solutions can be obtained by numerical methods in literature \cite{Bo}.

\section{Governing system of nonlinear elliptic equations}\label{s3}
\setcounter{equation}{0}

In this section, we derive the nonlinear elliptic equations to be studied and state our main result. Writing $g_{0}^{2}=g'^{2}=\xi=1$, we can render the nonlinear ordinary differential equations \eqref{2.18}--\eqref{2.21} into
\ber
\frac{f'}{r}-N\left[Q_{1}^{2}+(N-1)Q_{2}^{2}-N\right]&=&0,\label{3.1}\\
\frac{f'_{NA}}{r}-\frac{Q_{1}^{2}-Q_{2}^{2}}{2}&=&0,\label{3.2}\\
rQ'_{1}-Q_{1}\left(\frac{(N-1)f_{NA}+f}{N}\right)&=&0,\label{3.3}\\
rQ'_{2}-Q_{2}\left(\frac{-f_{NA}+f}{N}\right)&=&0.\label{3.4}
\eer
It is worth noting that $Q_{1}(r)\neq0$ for all $r\in(0,+\infty)$. If otherwise, $Q_{1}(r)=0$ at some $r\in(0,+\infty)$, then $Q_{1}(r)\equiv0$ for all $r\in(0,+\infty)$ by the continuous dependence theorem for solutions of the initial value problems of ordinary differential equations. Similarly, $Q_{2}(r)\neq0$ for all $r\in(0,+\infty)$. As a consequence, when $r>0$ the equations \eqref{3.3}--\eqref{3.4} are recast into
\ber
Nr\left(\mathrm{ln}Q_{1}\right)'&=&(N-1)f_{NA}+f,\label{3.7}\\
Nr\left(\mathrm{ln}Q_{2}\right)'&=&-f_{NA}+f.\label{3.8}
\eer
For convenience, we may introduce the new variable $u_{1}=\mathrm{ln}Q_{1}$, $u_{2}=\mathrm{ln}Q_{2}$. Then inserting the above equations into \eqref{3.1}--\eqref{3.2}, we arrive at the following equations
\ber
\frac{1}{r}\left(ru''_{1}+u'_{1}+(N-1)\left(ru''_{2}+u'_{2}\right)\right)&=&N\left(\e^{2u_{1}}+(N-1)\e^{2u_{2}}-N\right),\label{3.9}\\
\frac{1}{r}\left(ru''_{1}+u'_{1}-\left(ru''_{2}+u'_{2}\right)\right)&=&\frac{1}{2}\left(\e^{2u_{1}}-\e^{2u_{2}}\right),\label{3.10}
\eer
where $N>1$ is an integer. By a direct computation, we obtain the following system of nonlinear elliptic equations over $\mathbb{R}^{2}$
\ber
\triangle u_{1}(x)&=&\left(\frac{3}{2}-\frac{1}{2N}\right)\left(\mathrm{e}^{2u_{1}(x)}-1\right)+\left(N-\frac{3}{2}+\frac{1}{2N}\right)\left(\mathrm{e}^{2u_{2}(x)}-1\right)+4\pi n_{1}\delta(x),\label{3.11}\\
\triangle u_{2}(x)&=&\left(1-\frac{1}{2N}\right)\left(\mathrm{e}^{2u_{1}(x)}-1\right)+\left(N-1+\frac{1}{2N}\right)\left(\mathrm{e}^{2u_{2}(x)}-1\right)+4\pi n_{2}\delta(x),\label{3.12}
\eer
where $\delta(x)$ is the Dirac distribution on $\mathbb{R}^{2}$ concentrated at the origin and the positive integers $n_{1}$, $n_{2}$ are the multiplicities of the vortices corresponding to $u_{1}$ and $u_{2}$. We are interested in the existence of solutions of \eqref{3.11}--\eqref{3.12} for the case of the full plane. Therefore, we consider the system \eqref{3.11}--\eqref{3.12} over the plane with the topological boundary condition
\be
u_{1}(x)\rightarrow0,~~~u_{2}(x)\rightarrow0,~~~\text{as}~~~|x|\rightarrow\infty.\label{3.13}
\ee
Defining the matrix $A$
\be\label{3.14}
A=\left(
\begin{matrix}
\dfrac{3}{2}-\dfrac{1}{2N}&~~~N-\dfrac{3}{2}+\dfrac{1}{2N}\\[8pt]
1-\dfrac{1}{2N}&~~~
N-1+\dfrac{1}{2N}
\end{matrix}
\right)\triangleq
\left(
\begin{matrix}
\alpha&~~~\beta\\[8pt]
\alpha-\dfrac{1}{2}&~~~
\beta+\dfrac{1}{2}
\end{matrix}
\right),
\ee
then the equations \eqref{3.11}--\eqref{3.12} can be rewritten in a compact form
\be\label{3.15}
\triangle u_{i}=\sum_{j=1}^{2} a_{ij}\left(\mathrm{e}^{2u_{j}}-1\right)+4\pi n_{i}\delta(x),~~~i=1,2.
\ee
It is easily to see that our model has motion equations with the similar structure as those studied by Yang\cite{Yang}. We optimize the approach in \cite{Chen1,Chen2,Sakai,Yang,YangBook} to handle this system. In fact, letting $N=2$ in \eqref{3.15} gives the same equations as those arising in the generalized Abelian Higgs theory with $\left(U(1)^{m}\right)$ model studied in \cite{Yang,YangBook},
\ber
\triangle u_{1}&=&\frac{5}{4}\left(\mathrm{e}^{2u_{1}}-1\right)+\frac{3}{4}\left(\mathrm{e}^{2u_{2}}-1\right)+4\pi n_{1}\delta(x),\label{3.16}\\
\triangle u_{2}&=&\frac{3}{4}\left(\mathrm{e}^{2u_{1}}-1\right)+\frac{5}{4}\left(\mathrm{e}^{2u_{2}}-1\right)+4\pi n_{2}\delta(x).\label{3.17}
\eer

As shown in \cite{Yang}, it is possible to establish an existence and uniqueness theorem of the system \eqref{3.15} in $\mathbb{R}^{2}$ for more general matrices $A$. An explicit fulfillment of this idea is provided by our equations. Special attention should be paid here to the fact that even though the existence and uniqueness for solutions of the system \eqref{3.15} established in $\mathbb{R}^{2}$ is similar to the theorem shown in \cite{Yang}, our matrix $A$ does not satisfy the assumptions used to derive the decay estimates in \cite{Yang} since the matrix $A$ is not a positive definite real symmetric matrix. It should be emphasized that we can transform the real matrix $A$ into a positive definite real symmetric matrix thereby obtaining the demonstration of the decay estimates. Furthermore this decay estimates are closely related to the quantized integrals in the full plane. Next we will use the variational method of Jaffe and Taubes\cite{Jaffe} to get the existence for solutions of the nonlinear elliptic equations \eqref{3.11}--\eqref{3.12} over the full plane.

Concerning above situations, our main existence and uniqueness theorem for solutions of \eqref{3.11}--\eqref{3.12} are stated as follows.

\begin{theorem}\label{th3.1}
The nonlinear elliptic equations \eqref{3.11}--\eqref{3.12} over the full plane $\mathbb{R}^{2}$ subject to the topological boundary condition \eqref{3.13} always have a unique solution. Furthermore, this solution fulfills the boundary condition \eqref{3.13} exponentially fast. More precisely, for any small number $\varepsilon\in(0,1)$, there hold the following sharp decay estimates at infinity,
\ber
&&\left|pu_{1}(x)\right|^{2}+\left|2u_{2}(x)\right|^{2}\leq C(\varepsilon)\mathrm{e}^{-(1-\varepsilon)\sqrt{\lambda_{0}}|x|},\label{3.18}\\
&&\left|\nabla\left(mu_{1}(x)+2u_{2}(x)\right)\right|^{2}+\left|\nabla\left(pu_{1}(x)+qu_{2}(x)\right)\right|^{2}\leq C(\varepsilon)\mathrm{e}^{-(1-\varepsilon)\sqrt{\lambda}|x|},\label{3.19}
\eer
where
\ber\label{3.19a}
m=\frac{\left(2\alpha-1\right)^{2}}{2\beta\left(\lambda_{3}-\alpha\right)},~~p=\frac{2\alpha-1}{\beta},~~q=\frac{4\left(\lambda_{4}-\alpha\right)}{2\alpha-1},
\eer
$C(\varepsilon)$ is a positive constant depending only on $\varepsilon$ and $\lambda_{0},\lambda,\lambda_{3},\lambda_{4}$ are as defined by \eqref{4.39a}, \eqref{4.48} and \eqref{4.46b}. Besides, there hold the quantized integrals in the full plane,
\ber\label{3.20}
&&\int_{\mathbb{R}^{2}} \left\{\left[\left(m+2\right)\alpha-1\right]\left(\mathrm{e}^{2u_{1}(x)}-1\right)+\left[\left(m+2\right)\beta+1\right]\left(\mathrm{e}^{2u_{2}(x)}-1\right)\right\}\dd x=-4\pi\left(mn_{1}+2n_{2}\right),\nn\\[3mm]
&&\int_{\mathbb{R}^{2}} \left\{\left[\left(p+q\right)\alpha-\frac{q}{2}\right]\left(\mathrm{e}^{2u_{1}(x)}-1\right)+\left[\left(p+q\right)\beta+\frac{q}{2}\right]\left(\mathrm{e}^{2u_{2}(x)}-1\right)\right\}\dd x=-4\pi\left(pn_{1}+qn_{2}\right),\nn\\
\eer
where $\alpha$, $\beta$ are shown in \eqref{3.14} and $n_{1}$, $n_{2}$ are the multiplicities of the vortices corresponding to $u_{1}$ and $u_{2}$.
\end{theorem}

It is worth stressing that the facts stated in \eqref{3.20} actually appear in the form of the flux quantization in the corresponding quantum field theory model. The above theorem will be established in the subsequent sections.

\section{Proof of the main theorem}\label{s4}
\setcounter{equation}{0}

In this section, we study the existence and uniqueness of solutions to equations \eqref{3.11}--\eqref{3.12} under the topological boundary condition \eqref{3.13} and establish the decay estimates for the solutions, which allow us to obtain the quantized integrals over the full plane.

\subsection{Existence and uniqueness of the critical point}

As in \cite{Jaffe}, we need to take the background function
\be\label{4.0a}
u_{i}^{0}(x)=-n_{i}\mathrm{ln}\left(1+\tau|x|^{-2}\right),~~~i=1,2,~~~\tau>0,
\ee
then we can see that
\be\label{4.0b}
\triangle u_{i}^{0}(x)=-\varphi_{i}(x)+4\pi n_{i}\delta(x),~~~\varphi_{i}(x)=\frac{4n_{i}\tau}{\left(\tau+|x|^{2}\right)^{2}},~~~i=1,2.
\ee
It is significant noting that $\int_{\mathbb{R}^{2}}\varphi_{i}(x)\dd x=4\pi n_{i}(i=1,2)$. Setting $u_{i}(x)=u_{i}^{0}(x)+P_{i}(x)(i=1,2)$, the system \eqref{3.15} becomes
\ber
\triangle P_{1}(x)&\!=\!&\alpha\left(\mathrm{e}^{2(u_{1}^{0}(x)+P_{1}(x))}-1\right)\!+\!\beta\left(\mathrm{e}^{2(u_{2}^{0}(x)+P_{2}(x))}-1\right)+\varphi_{1}(x),\label{4.0c}\\[2mm]
\triangle P_{2}(x)&\!=\!&\left(\alpha-\frac{1}{2}\right)\left(\mathrm{e}^{2(u_{1}^{0}(x)+P_{1}(x))}-1\right)\!+\!\left(\beta+\frac{1}{2}\right)\left(\mathrm{e}^{2(u_{2}^{0}(x)+P_{2}(x))}-1\right)+\varphi_{2}(x).\label{4.0d}
\eer
To proceed, we use boldfaced letters to indicate column vectors in $\mathbb{R}^{2}$. With
\be\label{4.1}
\mathbf{P}=\left(P_{1},P_{2}\right)^{\tau},~~~\mathbf{E}=\left(\mathrm{e}^{2(u_{1}^{0}+P_{1})}-1,\mathrm{e}^{2(u_{2}^{0}+P_{2})}-1\right)^{\tau},~~~\mathbf{\Phi}=\left(\varphi_{1},\varphi_{2}\right)^{\tau}
\ee
and the notation set in \eqref{3.14}, the equations \eqref{4.0c}--\eqref{4.0d} can be written in the matrix form
\be\label{4.2}
\triangle\mathbf{P}=A\mathbf{E}+\mathbf{\Phi}.
\ee
This system is challenging since the coefficient matrix $A$ is not a positive definite or not even symmetric. Then in order to tackle it, we try to look for a variational principle.

To seek the variational principle, we will apply the property of the matrix $A$. It is clear that the matrix $A$ is simply nonsingular. According to the more general Crout decomposition theorem, we can see that there exist two $2\times2$ matrices, $L=\left(L_{jk}\right)$ which is lower triangular and $R=\left(R_{jk}\right)$ which is upper triangular, such that
\be\label{4.3}
A=LR.
\ee
Furthermore, by the scheme\cite{Stoer}, we can expictly constructed $L$ and $R$ from the coefficient matrix $A$ as follows
\ber
&&L=\left(
\begin{matrix}
1&~~~0\\[8pt]
\dfrac{2N-1}{3N-1}&~~~1
\end{matrix}
\right)=\left(
\begin{matrix}
1&~~~0\\[8pt]
1-\dfrac{1}{2\alpha}&~~~1
\end{matrix}
\right),\label{4.4}\\[4mm]
&&R=\left(
\begin{matrix}
\dfrac{3}{2}-\dfrac{1}{2N}&~~~
N-\dfrac{3}{2}+\dfrac{1}{2N}\\[8pt]
0&~~~\dfrac{N^{2}}{3N-1}
\end{matrix}
\right)=\left(
\begin{matrix}
\alpha&~~~
\beta\\[8pt]
0&~~~\dfrac{\alpha+\beta}{2\alpha}
\end{matrix}
\right),\label{4.5}
\eer
which will be conducive to the existence of a solution of the system \eqref{4.2}. With the new variable vector
\be\label{4.6}
\mathbf{w}=L^{-1}\mathbf{P}~~~\text{or}~~~\mathbf{P}=L\mathbf{w},
\ee
then the transformed system becomes
\be\label{4.7}
\triangle\mathbf{w}=R\mathbf{E}+L^{-1}\mathbf{\Phi}.
\ee
Moreover, it is convenient to rewrite the above system in the component form
\ber
\triangle w_{1}&=&\alpha\left(\mathrm{e}^{2\left(u_{1}^{0}+w_{1}\right)}-1\right)+\beta\left(\mathrm{e}^{2\left[u_{2}^{0}+\left(1-\frac{1}{2\alpha}\right)w_{1}+w_{2}\right]}-1\right)+\psi_{1},\label{4.8}\\
\triangle w_{2}&=&\frac{\alpha+\beta}{2\alpha}\left(\mathrm{e}^{2\left[u_{2}^{0}+\left(1-\frac{1}{2\alpha}\right)w_{1}+w_{2}\right]}-1\right)+\psi_{2},\label{4.9}
\eer
where $\psi_{1}=\varphi_{1}$, $\psi_{2}=(\frac{1}{2\alpha}-1)\varphi_{1}+\varphi_{2}$. In order to see the variational structure of \eqref{4.8}--\eqref{4.9} clearly, it is beneficial to reformulate them equivalently as
\ber
\frac{2\alpha-1}{\alpha\beta}\triangle w_{1}&=&\frac{2\alpha-1}{\beta}\left(\mathrm{e}^{2\left(u_{1}^{0}+w_{1}\right)}-1\right)+\left(2-\frac{1}{\alpha}\right)\left(\mathrm{e}^{2\left[u_{2}^{0}+\left(1-\frac{1}{2\alpha}\right)w_{1}+w_{2}\right]}-1\right)\nn\\
&+&\frac{2\alpha-1}{\alpha\beta}\psi_{1},\label{4.10}\\
\frac{4\alpha}{\alpha+\beta}\triangle w_{2}&=&2\left(\mathrm{e}^{2\left[u_{2}^{0}+\left(1-\frac{1}{2\alpha}\right)w_{1}+w_{2}\right]}-1\right)+\frac{4\alpha}{\alpha+\beta}\psi_{2}.\label{4.11}
\eer

To accommodate the boundary condition \eqref{3.13}, we will work on the standard Sobolev space $W^{1,2}(\mathbb{R}^2)\times W^{1,2}(\mathbb{R}^2)$. It is obvious to see that equations \eqref{4.10}--\eqref{4.11} are the Euler-Lagrange equations of the following action functional
\ber\label{4.12}
&&I(w_{1},w_{2})=\int_{\mathbb{R}^{2}}\Bigg\{\frac{2\alpha-1}{2\alpha\beta}\left|\nabla w_{1}\right|^{2}+\frac{2\alpha}{\alpha+\beta}\left|\nabla w_{2}\right|^{2}+\mathrm{e}^{2u_{2}^{0}}\left(\mathrm{e}^{2\left[\left(1-\frac{1}{2\alpha}\right)w_{1}+w_{2}\right]}-1\right)\nn\\
&&~~~~~~~~~~~~~~~~~~~~~~~+\frac{2\alpha-1}{2\beta}\mathrm{e}^{2u_{1}^{0}}\left(\mathrm{e}^{2w_{1}}-1\right)+\frac{2\alpha-1}{\alpha\beta}\psi_{1}w_{1}-\frac{\left(2\alpha-1\right)\left(\alpha+\beta\right)}{\alpha\beta}w_{1}\nn\\
&&~~~~~~~~~~~~~~~~~~~~~~~+\frac{4\alpha}{\alpha+\beta}\psi_{2}w_{2}-2w_{2}\Bigg\}\dd x.
\eer
There is no difficulty in checking that the functional $I$ is indeed a $C^{1}$-functional for $w_{1}$, $w_{2}\in W^{1,2}(\mathbb{R}^2)$. Then we need only to find the critical points of the functional $I$ defined in \eqref{4.12} for purpose of solving the equations \eqref{4.8}--\eqref{4.9}. To this end, we apply a direct method developed in \cite{Chen2}.

To obtain the critical points of the functional $I$, our first step is to show that it is coercive over $W^{1,2}(\mathbb{R}^2)$. The form of \eqref{4.12} allows us to get that its Fr$\acute{e}$chet derivative satisfies
\ber\label{4.13}
&&DI(w_{1},w_{2})(w_{1},w_{2})
=\int_{\mathbb{R}^{2}}\Bigg\{\frac{2\alpha-1}{\alpha\beta}\left|\nabla w_{1}\right|^{2}+\frac{4\alpha}{\alpha+\beta}\left|\nabla w_{2}\right|^{2}+\left(\mathrm{e}^{2u_{2}^{0}}-1\right)\left[\left(2-\frac{1}{\alpha}\right)w_{1}+2w_{2}\right]\nn\\
&&~~~~~~~~~~~~~~~~~~~~~~~~~~~~~~~~~~~~+\frac{2\alpha-1}{\beta}\left(\mathrm{e}^{2u_{1}^{0}}-1\right)w_{1}+\frac{2\alpha-1}{\beta}\mathrm{e}^{2u_{1}^{0}}\left(\mathrm{e}^{2w_{1}}-1\right)w_{1}\nn\\
&&~~~~~~~~~~~~~~~~~~~~~~~~~~~~~~~~~~~~+\mathrm{e}^{2u_{2}^{0}}\left(\mathrm{e}^{\left(2-\frac{1}{\alpha}\right)w_{1}+2w_{2}}-1\right)\left[\left(2-\frac{1}{\alpha}\right)w_{1}+2w_{2}\right]\nn\\
&&~~~~~~~~~~~~~~~~~~~~~~~~~~~~~~~~~~~~+\frac{2\alpha-1}{\alpha\beta}\psi_{1}w_{1}+\frac{4\alpha}{\alpha+\beta}\psi_{2}w_{2}\Bigg\}\dd x.
\eer

Noting that
\ber\label{4.14}
&&\left|\nabla w_{1}\right|^{2}+\left|\nabla w_{2}\right|^{2}=\left|\nabla P_{1}\right|^{2}+\left|\nabla \left(P_{2}-\gamma P_{1}\right)\right|^{2}\nn\\
&&~~~~~~~~~~~~~~~~~~~~~\,\leq\left|\nabla P_{1}\right|^{2}+\left|\nabla P_{2}\right|^{2}+\gamma^{2}\left|\nabla P_{1}\right|^{2}+2\gamma\left|\left(\nabla P_{1}, \nabla P_{2}\right)\right|\nn\\
&&~~~~~~~~~~~~~~~~~~~~~\,\leq\left(1+\gamma+\gamma^{2}\right)\left|\nabla P_{1}\right|^{2}+\left(1+\gamma\right)\left|\nabla P_{2}\right|^{2}\nn\\
&&~~~~~~~~~~~~~~~~~~~~~\,<3\left(\left|\nabla P_{1}\right|^{2}+\left|\nabla P_{2}\right|^{2}\right),
\eer
where $\gamma=\frac{2N-1}{3N-1}\in\left(\frac{1}{2},\frac{2}{3}\right)$. On the other hand, we have
\ber\label{4.15}
&&\left|\nabla w_{1}\right|^{2}+\left|\nabla w_{2}\right|^{2}\geq\left(1+\gamma^{2}\right)\left|\nabla P_{1}\right|^{2}+\left|\nabla P_{2}\right|^{2}-2\gamma\left|\left(\nabla P_{1}, \nabla P_{2}\right)\right|\nn\\
&&~~~~~~~~~~~~~~~~~~~~~\,\geq\left(1+\gamma^{2}-\frac{\gamma}{2\varepsilon}\right)\left|\nabla P_{1}\right|^{2}+\left(1-2\gamma\varepsilon\right)\left|\nabla P_{2}\right|^{2}\nn\\
&&~~~~~~~~~~~~~~~~~~~~~\,\geq\left(\frac{5}{4}-\frac{1}{3\varepsilon}\right)\left|\nabla P_{1}\right|^{2}+\left(1-\frac{4\varepsilon}{3}\right)\left|\nabla P_{2}\right|^{2},
\eer
for any $\varepsilon\in\left(\frac{4}{15},\frac{3}{4}\right)$. Taking $\varepsilon=\frac{1}{3}$, then
\be\label{4.16}
\left|\nabla w_{1}\right|^{2}+\left|\nabla w_{2}\right|^{2}\geq\frac{1}{4}\left|\nabla P_{1}\right|^{2}+\frac{5}{9}\left|\nabla P_{2}\right|^{2}>\frac{5}{9}\left(\left|\nabla P_{1}\right|^{2}+\left|\nabla P_{2}\right|^{2}\right).
\ee
Combining \eqref{4.14} and \eqref{4.16}, we can get
\be\label{4.17}
\frac{5}{9}\left(\left|\nabla P_{1}\right|^{2}+\left|\nabla P_{2}\right|^{2}\right)<\left|\nabla w_{1}\right|^{2}+\left|\nabla w_{2}\right|^{2}<3\left(\left|\nabla P_{1}\right|^{2}+\left|\nabla P_{2}\right|^{2}\right).
\ee
Similarly, it can be inferred that
\be\label{4.18}
\frac{5}{9}\left(P_{1}^{2}+P_{2}^{2}\right)<w_{1}^{2}+w_{2}^{2}<3\left(P_{1}^{2}+P_{2}^{2}\right).
\ee
In view of \eqref{4.13}, \eqref{4.17} and \eqref{4.18}, we observe that
\ber\label{4.19}
&&DI(w_{1},w_{2})(w_{1},w_{2})-\frac{5}{9}\min\left\{\frac{2\alpha-1}{2\beta},\frac{4\alpha}{\alpha+\beta}\right\}\int_{\mathbb{R}^{2}}\left(\left|\nabla P_{1}\right|^{2}+\left|\nabla P_{2}\right|^{2}\right)\dd x\nn\\[2mm]
&&~~~~~~~~~~~~~~~~~~~~~~~~~\geq\int_{\mathbb{R}^{2}}\Bigg\{\frac{2\alpha-1}{\beta}\left(\mathrm{e}^{2(u_{1}^{0}+P_{1})}-1\right)P_{1}+2\left[\frac{2\alpha-1}{2\alpha\beta}\psi_{1}+\frac{1-2\alpha}{\alpha+\beta}\psi_{2}\right]P_{1}\nn\\
&&~~~~~~~~~~~~~~~~~~~~~~~~~~~~~~~~~~~~+2\left(\mathrm{e}^{2(u_{2}^{0}+P_{2})}-1\right)P_{2}+\frac{4\alpha}{\alpha+\beta}\psi_{2}P_{2}\Bigg\}\dd x\nn\\[2mm]
&&~~~~~~~~~~~~~~~~~~~~~~~~~=\frac{2\alpha-1}{2\beta}\int_{\mathbb{R}^{2}}2P_{1}\left(\mathrm{e}^{2(u_{1}^{0}+P_{1})}-1+X_{1}\right)\dd x+\int_{\mathbb{R}^{2}}2P_{2}\left(\mathrm{e}^{2(u_{2}^{0}+P_{2})}-1+X_{2}\right)\dd x\nn\\[2mm]
&&~~~~~~~~~~~~~~~~~~~~~~~~~\equiv \frac{2\alpha-1}{2\beta}M_{1}(P_{1})+M_{2}(P_{2}),
\eer
where
\be
X_{1}=\frac{1}{\alpha}\psi_{1}-\frac{2\beta}{\alpha+\beta}\psi_{2},~~~X_{2}=\frac{2\alpha}{\alpha+\beta}\psi_{2}. \nn\\
\ee
Next we need to estimate the term $M_{i}(P_{i})(i=1,2)$ on the right-hand side of the above. To proceed further, as in \cite{Jaffe}, we can choose a decomposition $P_{i}=P_{i+}-P_{i-}(i=1,2)$ with $P_{+}=\max\{0,P\}$ and $P_{-}=\max\{0,-P\}$ for $P\in\mathbb{R}$. Then $M_{i}(P_{i})=M_{i}(P_{i+})+M_{i}(-P_{i-})$.

According to the elementary inequality $\mathrm{e}^{t}-1\geq t$ for $t\in\mathbb{R}$, we can obtain
\ber\label{4.20}
&&M_{1}(P_{1+})\geq\int_{\mathbb{R}^{2}}(2P_{1+})^{2}\dd x+\int_{\mathbb{R}^{2}}2P_{1+}\left(2u_{1}^{0}+X_{1}\right)\dd x\nn\\[2mm]
&&~~~~~~~~~~~\geq\frac{1}{2}\left\|2P_{1+}\right\|_{2}^{2}-\frac{1}{2}\int_{\mathbb{R}^{2}}\left(2u_{1}^{0}+X_{1}\right)^{2}\dd x\nn\\[2mm]
&&~~~~~~~~~~~\geq\frac{1}{2}\int_{\mathbb{R}^{2}}\frac{\left(2P_{1+}\right)^{2}}{1+2P_{1+}}\dd x-C,
\eer
where we have used the fact $u_{1}^{0}, X_{1}\in L^{2}(\mathbb{R}^{2})$. Here and what follows we use $C$ to denote a general positive constant that can take different values at different places. Then applying the inequality
\be
1-\mathrm{e}^{-t}\geq\frac{t}{1+t},~~~ \forall t\geq0,\nn\\
\ee
we can estimate $M_{1}(-P_{1-})$ as follows
\ber\label{4.21}
&&M_{1}(-P_{1-})=\int_{\mathbb{R}^{2}}2P_{1-}\left(1-\mathrm{e}^{2(u_{1}^{0}-P_{1-})}-X_{1}\right)\dd x\nn\\[2mm]
&&~~~~~~~~~~~~\,~=\int_{\mathbb{R}^{2}}2P_{1-}\left[\mathrm{e}^{2u_{1}^{0}}\left(1-\mathrm{e}^{-2P_{1-}}\right)-\mathrm{e}^{2u_{1}^{0}}+1-X_{1}\right]\dd x\nn\\[2mm]
&&~~~~~~~~~~~~\,~\geq\int_{\mathbb{R}^{2}}2P_{1-}\left(\mathrm{e}^{2u_{1}^{0}}\frac{2P_{1-}}{1+2P_{1-}}-\mathrm{e}^{2u_{1}^{0}}+1-X_{1}\right)\dd x\nn\\[2mm]
&&~~~~~~~~~~~~\,~=\int_{\mathbb{R}^{2}}\frac{2P_{1-}}{1+2P_{1-}}\left[\left(1+2P_{1-}\right)\left(1-X_{1}-\mathrm{e}^{2u_{1}^{0}}\right)+2P_{1-}\mathrm{e}^{2u_{1}^{0}}\right]\dd x\nn\\[2mm]
&&~~~~~~~~~~~~\,~=\int_{\mathbb{R}^{2}}\frac{2P_{1-}}{1+2P_{1-}}\left(1-X_{1}-\mathrm{e}^{2u_{1}^{0}}\right)\dd x+\int_{\mathbb{R}^{2}}\frac{\left(2P_{1-}\right)^{2}}{1+2P_{1-}}\left(1-X_{1}\right)\dd x.
\eer
Noting the definition of \eqref{4.0b}, then we can choose $\tau>0$ sufficiently large such that
\be
X_{1}(x)<\frac{1}{2},~~~\forall x\in\mathbb{R}^{2}.\nn\\
\ee
In view of the fact $1-\mathrm{e}^{2u_{1}^{0}}$ and $X_{1}$ both belong to $L^{2}(\mathbb{R}^{2})$, we see that
\ber\label{4.21a}
\int_{\mathbb{R}^{2}}\frac{2P_{1-}}{1+2P_{1-}}\left|1-X_{1}-\mathrm{e}^{2u_{1}^{0}}\right|\dd x\leq\varepsilon\int_{\mathbb{R}^{2}}\frac{\left(2P_{1-}\right)^{2}}{(1+2P_{1-})^{2}}\dd x+C(\varepsilon),
\eer
where $\varepsilon>0$ could be taken to be arbitrarily small. Therefore, combining \eqref{4.21a} and \eqref{4.21}, we find that
\be\label{4.21b}
M_{1}(-P_{1-})\geq\frac{1}{4}\int_{\mathbb{R}^{2}}\frac{\left(2P_{1-}\right)^{2}}{(1+2P_{1-})^{2}}\dd x-C_{1}(\varepsilon),
\ee
provided that $\varepsilon<\frac{1}{4}$. Recall the lower estimate for $M_{1}(P_{1+})$ obtained earlier. As a consequence, we find that
\be\label{4.22}
M_{1}(P_{1})\geq\frac{1}{4}\int_{\mathbb{R}^{2}}\frac{\left(2P_{1}\right)^{2}}{(1+|2P_{1}|)^{2}}\dd x-C.
\ee
Analogous estimates could be made for  $M_{2}(P_{2})$
\be\label{4.23}
M_{2}(P_{2})\geq\frac{1}{4}\int_{\mathbb{R}^{2}}\frac{\left(2P_{2}\right)^{2}}{(1+|2P_{2}|)^{2}}\dd x-C.
\ee
Hence, from \eqref{4.22} and \eqref{4.23}, we arrive at
\ber\label{4.24}
&&DI(w_{1},w_{2})(w_{1},w_{2})-\frac{5}{9}\min\left\{\frac{2\alpha-1}{2\beta},\frac{4\alpha}{\alpha+\beta}\right\}\int_{\mathbb{R}^{2}}\left(\left|\nabla P_{1}\right|^{2}+\left|\nabla P_{2}\right|^{2}\right)\dd x\nn\\
&&~~~~~~~~~~~~~~~~~~~~~~~~~\geq\frac{2\alpha-1}{8\beta}\int_{\mathbb{R}^{2}}\frac{\left(2P_{1}\right)^{2}}{(1+|2P_{1}|)^{2}}\dd x+\frac{1}{4}\int_{\mathbb{R}^{2}}\frac{\left(2P_{2}\right)^{2}}{(1+|2P_{2}|)^{2}}\dd x-C.
\eer

In order to analyze $\left\|2P_{i}\right\|_{L^{2}(\mathbb{R}^{2})}(i=1,2)$, we now stress the Gagliardo--Nirenberg--Sobolev embedding inequality\cite{Gilbarg} over $W^{1,2}(\mathbb{R}^2)$
\be\label{4.25}
\int_{\mathbb{R}^{2}}f^{4}\dd x\leq2\int_{\mathbb{R}^{2}}f^{2}\dd x\int_{\mathbb{R}^{2}}|\nabla f|^{2}\dd x,~~\forall f\in W^{1,2}(\mathbb{R}^2).
\ee
We will use \eqref{4.25} to prove the desired coercivity inequality. In reality, by virtue of \eqref{4.25}, we can obtain
\ber\label{4.26}
&&\left(\int_{\mathbb{R}^{2}}(2P_{i})^{2}\dd x\right)^{2}=\left(\int_{\mathbb{R}^{2}}\frac{|2P_{i}|}{1+|2P_{i}|}\left(1+|2P_{i}|\right)|2P_{i}|\dd x\right)^{2}\nn\\
&&~~~~~~~~~~~~~~~~~~~~~\leq\int_{\mathbb{R}^{2}}\frac{\left(2P_{i}\right)^{2}}{\left(1+|2P_{i}|\right)^{2}}\dd x\int_{\mathbb{R}^{2}}\left(1+|2P_{i}|\right)^{2}|2P_{i}|^{2}\dd x\nn\\
&&~~~~~~~~~~~~~~~~~~~~~\leq4\int_{\mathbb{R}^{2}}\frac{\left(2P_{i}\right)^{2}}{\left(1+|2P_{i}|\right)^{2}}\dd x\int_{\mathbb{R}^{2}}|2P_{i}|^{2}\dd x\left(1+\int_{\mathbb{R}^{2}}|2\nabla P_{i}|^{2}\dd x\right)\nn\\
&&~~~~~~~~~~~~~~~~~~~~~\leq\frac{1}{2}\left(\int_{\mathbb{R}^{2}}|2P_{i}|^{2}\dd x\right)^{2}+32\Bigg\{1+\left[\int_{\mathbb{R}^{2}}\frac{\left(2P_{i}\right)^{2}}{\left(1+|2P_{i}|\right)^{2}}\dd x\right]^{4}\nn\\
&&~~~~~~~~~~~~~~~~~~~~~~~~+\left(\int_{\mathbb{R}^{2}}|2\nabla P_{i}|^{2}\dd x\right)^{4}\Bigg\},~~i=1,2
\eer
which yields
\ber\label{4.27}
\left\|2P_{i}\right\|_{L^{2}(\mathbb{R}^{2})}\leq C\left[1+\int_{\mathbb{R}^{2}}\frac{\left(2P_{i}\right)^{2}}{\left(1+|2P_{i}|\right)^{2}}\dd x+\int_{\mathbb{R}^{2}}|2\nabla P_{i}|^{2}\dd x\right],~~i=1,2,
\eer
where $C$ denote a positive constant. Therefore, inserting \eqref{4.17}, \eqref{4.18} and \eqref{4.27} into \eqref{4.24}, we get
\be\label{4.28}
DI(w_{1},w_{2})(w_{1},w_{2})\geq C_{1}\left(\left\|w_{1}\right\|_{W^{1,2}(\mathbb{R}^{2})}+\left\|w_{2}\right\|_{W^{1,2}(\mathbb{R}^{2})}\right)-C_{2}
\ee
for suitable positive constants $C_{1}$ and $C_{2}$, which gives the expected coerciveness of the functional $I$ over $W^{1,2}(\mathbb{R}^2)$. As a result of the coercive lower bound \eqref{4.28}, we can now infer that the action functional $I$ admits a critical point in the space $W^{1,2}(\mathbb{R}^2)$, which follows in a standard path\cite{YangBook}.

Our next step is to utilize \eqref{4.28} to show that the system \eqref{4.8}--\eqref{4.9} has a solution by confirming that \eqref{4.12} has a critical point. In fact, in view of \eqref{4.28}, we can let $R>0$ large enough such that
\be\label{4.29}
\inf\left\{DI(w_{1},w_{2})\Big|w_{1},w_{2}\in W^{1,2}(\mathbb{R}^2),\left\|w_{1}\right\|_{W^{1,2}(\mathbb{R}^2)}+\left\|w_{2}\right\|_{W^{1,2}(\mathbb{R}^2)}=R\right\}\geq1
\ee
and consider the optimization problem
\be\label{4.30}
\eta=\min\left\{I(w_{1},w_{2})\Big|\left\|w_{1}\right\|_{W^{1,2}(\mathbb{R}^2)}+\left\|w_{2}\right\|_{W^{1,2}(\mathbb{R}^2)}\leq R\right\}.
\ee
Let $\{(w_{1}^{(n)},w_{2}^{(n)})\}$ be a minimizing sequence of the problem \eqref{4.30}. Since $\{(w_{1}^{(n)},w_{2}^{(n)})\}$ is bounded in $W^{1,2}(\mathbb{R}^2)$, then its subsequence is weakly convergent. Without loss of generality, we may assume $\{(w_{1}^{(n)},w_{2}^{(n)})\}$ is also weakly converges to $\{(w_{1},w_{2})\}$ in $W^{1,2}(\mathbb{R}^2)$. It is worth noting that the functional \eqref{4.12} is continuous, differentiable and convex in $W^{1,2}(\mathbb{R}^2)$, thus the functional $I$ is weakly lower semi-continuous. According to the Fatou lemma, we have $I(w_{1},w_{2})\leq\lim\limits_{n\rightarrow\infty}I(w_{1}^{(n)},w_{2}^{(n)})=\eta$. On the other hand, since the norm of $W^{1,2}(\mathbb{R}^2)$ is also weakly lower semicontinuous, we obtain $\|w_{1}\|_{W^{1,2}(\mathbb{R}^2)}+\|w_{2}\|_{W^{1,2}(\mathbb{R}^2)}\leq R$, that is $I(w_{1},w_{2})\in\eta$. Hence $I(w_{1},w_{2})\geq\eta$. In summary, $I(w_{1},w_{2})=\eta$. In other words, the minimization problem \eqref{4.30} admits a solution $(w_{1},w_{2})$.

In the following, to show that $(w_{1},w_{2})$ is a critical point of $I$ which is a weak solution of equations \eqref{4.8}--\eqref{4.9}, we only need to argue that $(w_{1},w_{2})$ is an interior point or
\be
\|w_{1}\|_{W^{1,2}(\mathbb{R}^2)}+\|w_{2}\|_{W^{1,2}(\mathbb{R}^2)}<R.\nn
\ee
Suppose by contradiction that $\|w_{1}\|_{W^{1,2}(\mathbb{R}^2)}+\|w_{2}\|_{W^{1,2}(\mathbb{R}^2)}=R$. Since
\be\label{4.31}
\left\|\left(w_{1},w_{2}\right)-t\left(w_{1},w_{2}\right)\right\|_{W^{1,2}(\mathbb{R}^2)}=(1-t)R<R,~~\forall t\in(0,1),
\ee
that is, $(w_{1}^{t},w_{2}^{t})=(1-t)(w_{1},w_{2})$ is an interior point for any $t\in(0,1)$, then we can get
\be\label{4.32}
I\left(w_{1}^{t},w_{2}^{t}\right)\geq I\left(w_{1},w_{2}\right)=\eta.
\ee
However, due to \eqref{4.29}, we arrive that
\be\label{4.33}
\lim\limits_{t\rightarrow0}\frac{I\left(w_{1}^{t},w_{2}^{t}\right)-I\left(w_{1},w_{2}\right)}{t}=\frac{\dd}{\dd t}\left(I\left(w_{1}^{t},w_{2}^{t}\right)\right)\Big|_{t=0}=-DI\left(w_{1},w_{2}\right)\left(w_{1},w_{2}\right)\leq-1.
\ee
Consequently, when $t>0$ is sufficiently small, in virtue of \eqref{4.33}, we can know that
\be\label{4.34}
I\left(w_{1}^{t},w_{2}^{t}\right)<I\left(w_{1},w_{2}\right)=\eta,
\ee
which contradicts \eqref{4.32}. Thus we can attain that $(w_{1},w_{2})$ must be an interior point for the problem \eqref{4.30} which as a critical point of $I$ in $W^{1,2}({\mathbb{R}^{2}})$ solves the system \eqref{4.8}--\eqref{4.9}.

Finally, the strict convexity of $I$ already says that such a critical point must be unique. Therefore, the functional $I$ can only have at most one critical point in the space $W^{1,2}({\mathbb{R}^{2}})$ and uniqueness of a solution to \eqref{4.8}--\eqref{4.9} follows. Indeed, such a uniqueness outcome follows in a more straightforward manner from the structure of the functional. It is worth to emphasize that the part in the integrand of the functional \eqref{4.12} can be recast in addition to the derivative terms of $(w_{1},w_{2})$ as
\ber
&&F\left(w_{1},w_{2}\right)=\mathrm{e}^{2u_{2}^{0}}\left(\mathrm{e}^{2\left[\left(1-\frac{1}{2\alpha}\right)w_{1}+w_{2}\right]}-1\right)+\frac{2\alpha-1}{2\beta}\mathrm{e}^{2u_{1}^{0}}\left(\mathrm{e}^{2w_{1}}-1\right)+\frac{2\alpha-1}{\alpha\beta}\psi_{1}w_{1}\nn\\
&&~~~~~~~~~~~~~~~~~-\frac{\left(2\alpha-1\right)\left(\alpha+\beta\right)}{\alpha\beta}w_{1}+\frac{4\alpha}{\alpha+\beta}\psi_{2}w_{2}-2w_{2},
\eer
whose Hessian matrix is not hard verified to be positive definite. As a consequence, the functional $I$ is strictly convex, which implies $(w_{1},w_{2})$ is the unique critical point of $I$ in the space $W^{1,2}({\mathbb{R}^{2}})$. Then the existence and uniqueness of a critical point of $I$ in $W^{1,2}({\mathbb{R}^{2}})$ is obtained. Furthermore, this critical point is a smooth solution of the system \eqref{4.8}--\eqref{4.9} in view of the standard elliptic regularity theory.

\subsection{Asymptotic behavior and quantized integrals}

Now we research the behavior the solution at infinity. For the solution gained in the previous part, we first derive some pointwise decay properties by elliptic $L^{p}$--estimates and the standard embedding inequalities, and then estimate the asymptotic decay rates of the solutions and their derivatives near infinity through the maximum principle. As an application of the decay estimate, we can compute the quantized integrals described in Theorem \ref{th3.1}.

In the following, we claim that if $h\in W^{1,2}({\mathbb{R}^{2}})$, then $\mathrm{e}^{h}-1\in L^{2}(\mathbb{R}^{2})$. Since the Sobolev embedding inequality in two dimensions
\be\label{4.35}
\left\|h\right\|_{L^{k}(\mathbb{R}^{2})}\leq\left(\pi\left(\frac{k-2}{2}\right)\right)^{\frac{k-2}{2k}}\left\|h\right\|_{W^{1,2}(\mathbb{R}^{2})},~~h\in W^{1,2}(\mathbb{R}^{2}),~~k>2
\ee
and the MacLaurin series
\be\label{4.36}
\left(\mathrm{e}^{h}-1\right)^{2}=h^{2}+\sum_{k=3}^{\infty} \frac{2^{k}-2}{k!}h^{k},
\ee
we have seen that
\be\label{4.37}
\left\|\mathrm{e}^{h}-1\right\|_{L^{2}(\mathbb{R}^{2})}^{2}\leq\left\|h\right\|_{L^{k}(\mathbb{R}^{2})}^{2}+\sum_{k=3}^{\infty} \frac{2^{k}-2}{k!}\left(\pi\left(\frac{k-2}{2}\right)\right)^{\frac{k-2}{2}}\left\|h\right\|_{W^{1,2}(\mathbb{R}^{2})}^{k}.
\ee
Assume $\alpha_{k}=\frac{2^{k}-2}{k!}\left(\pi\left(\frac{k-2}{2}\right)\right)^{\frac{k-2}{2}}\left\|h\right\|_{W^{1,2}(\mathbb{R}^{2})}^{k}$, by virtue of the Stirling formula
\be
n!\thicksim\sqrt{2\pi n}\left(\frac{n}{\mathrm{e}}\right)^{n}~~(n\rightarrow\infty),\nn\\
\ee
then by straightforward calculations we get
\ber
&&\sqrt[k]{\alpha_{k}}\thicksim\frac{\sqrt[k]{2^{k}-2}}{\mathrm{e}^{-1}k\left(2k\pi\right)^{\frac{1}{2k}}}\left(\frac{k-2}{2}\pi\right)^{\frac{k-2}{2k}}\left\|h\right\|_{W^{1,2}(\mathbb{R}^{2})}\nn\\
&&~~~~~~\thicksim 2\mathrm{e}\sqrt{\pi}\left(\frac{k-2}{2k^{2}}\right)^{\frac{1}{2}}\left\|h\right\|_{W^{1,2}(\mathbb{R}^{2})}\rightarrow0~~(k\rightarrow\infty).\nn
\eer
Whence, it is easy to observe that the series on the right--hand side of \eqref{4.37} is convergent which proves our claim.

It is worth noting that $w_{1},w_{2}\in W^{1,2}(\mathbb{R}^{2})$, applying the above claim we can infer that the right--hand side of \eqref{4.8}--\eqref{4.9} all lie in $L^{2}(\mathbb{R}^{2})$, which establishes $w_{1},w_{2}\in W^{2,2}(\mathbb{R}^{2})$ by the well-known $L^{2}$--estimate for elliptic equations. In addition, according to the standard Sobolev embeddings and the fact that we are in two dimensions, we observe that $w_{1}(x),w_{2}(x)\rightarrow0$ as $|x|\rightarrow\infty$ which gives the desired boundary condition \eqref{3.13} at infinity. Due to the fact $w_{1},w_{2}\in W^{2,2}(\mathbb{R}^{2})$ and the embedding $W^{1,2}(\mathbb{R}^{2})\hookrightarrow L^{p}(\mathbb{R}^{2})(p>2)$, we can obtain that the right--hand side of \eqref{4.8}--\eqref{4.9} all belongs to $L^{p}(\mathbb{R}^{2})$ for any $p>2$. Furthermore, the elliptic $L^{p}$ estimate enable us to get that $w_{1},w_{2}\in W^{2,p}(\mathbb{R}^{2})(p>2)$. Consequently, $\left|\nabla w_{1}\right|(x)\rightarrow0,\left|\nabla w_{2}\right|(x)\rightarrow0$ when $|x|\rightarrow\infty$, as expected.

Finally, we estimate the decay rates for $u_{1},u_{2}$ and $|\nabla(mu_{1}+2u_{2})|, |\nabla(pu_{1}+qu_{2})|$, where $m, p, q$ are as defined by \eqref{3.19a}. For the coefficient matrix $A$ of the system \eqref{3.15}, there exists a diagonal matrix $B$ such that the matrix $BA=M$ is positive definite and symmetric, where
\ber
&&B=\left(b_{ij}\right)=\left(
\begin{matrix}
\dfrac{2\alpha-1}{\beta}&~~~0\\[8pt]
0&~~~
2
\end{matrix}
\right),
B^{-1}=\left(b_{ij}^{-1}\right)=\left(\frac{1}{b_{ij}}\right)=\left(
\begin{matrix}
\dfrac{\beta}{2\alpha-1}&~~~0\\[8pt]
0&~~~
\dfrac{1}{2}
\end{matrix}
\right),~~i,j=1,2,\nn\\
&&BA=\left(
\begin{matrix}
\dfrac{2\alpha-1}{\beta}&~~~0\\[8pt]
0&~~~
2
\end{matrix}
\right)
\left(
\begin{matrix}
\alpha&~~~\beta\\[8pt]
\alpha-\dfrac{1}{2}&~~~
\beta+\dfrac{1}{2}
\end{matrix}
\right)
=\left(
\begin{matrix}
\dfrac{2\alpha^{2}-\alpha}{\beta}&~~~2\alpha-1\\[8pt]
2\alpha-1&~~~
2\beta+1
\end{matrix}
\right)\triangleq M.\nn
\eer
Therefore, we reduce \eqref{3.11}--\eqref{3.12} in $\mathbb{R}^{2}\backslash(0,0)$ to
\be\label{4.38}
\triangle v_{j}=\sum_{k=1}^{2} m_{jk}v_{k}+\sum_{k=1}^{2} m_{jk}\left(\mathrm{e}^{2\sum\limits_{l=1}^{2} b_{kl}^{-1}v_{k}}-1-v_{k}\right),~~~j=1,2.
\ee

Setting $O$ be a $2\times2$ orthogonal matrix such that
\be\label{4.39}
O^{\tau}MO=\mathrm{diag}\left\{\lambda_{1},\lambda_{2}\right\},~~\lambda_{1},\lambda_{2}>0,
\ee
where
\be\label{4.39a}
\lambda_{0}=\min\left\{\lambda_{1},\lambda_{2}\right\}>0.
\ee
Hence, as to the new variable vector
\be\label{4.40}
\mathbf{g}=\left(g_{1},g_{2}\right)^{\tau}=O^{\tau}\left(v_{1},v_{2}\right)^{\tau}=O^{\tau}\mathbf{v},
\ee
substitute \eqref{4.40} into \eqref{4.38}, then according to the \eqref{4.39} and the behaviour of $\mathbf{g}(x)\rightarrow0$ as $|x|\rightarrow\infty$, we have
\be\label{4.41}
\triangle g_{j}=\lambda_{j}g_{j}+\sum_{k=1}^{2} C_{jk}(x)g_{k},~~j=1,2,
\ee
where $C_{jk}(x)(j,k=1,2)$ depend on $\mathbf{g}(x)$ and $C_{jk}(x)\rightarrow0$ as $|x|\rightarrow\infty(j,k=1,2)$. Let $g^{2}=g_{1}^{2}+g_{2}^{2}$, then we arrive that
\be\label{4.42}
\triangle g^{2}\geq \lambda_{0}g^{2}-C(x)g^{2},~~x\in\mathbb{R}^{2}\backslash(0,0),
\ee
where $C(x)\rightarrow 0$ as $|x|\rightarrow\infty$. Thus, for any sufficiently small $\varepsilon\in(0,1)$, there is a suitably large $R_{\varepsilon}>R$ so that
\be\label{4.43}
\triangle g^{2}\geq \left(1-\frac{\varepsilon}{2}\right)\lambda_{0}g^{2},~~x\in\mathbb{R}^{2}\backslash(0,0).
\ee
As a consequence, in view of a comparison function argument and the property $g_{1}^{2}+g_{2}^{2}=0$ at infinity, we attain a positive constant $C(\varepsilon)$ to make
\be\label{4.44}
|\mathbf{v}|^{2}=|\mathbf{g}|^{2}=g^{2}\leq C(\varepsilon)\mathrm{e}^{-(1-\varepsilon)\sqrt{\lambda_{0}}|x|},~~\text{as}~~|x|\to \infty
\ee
valid, which leads to the desired decay estimates \eqref{3.18} stated in Theorem \ref{th3.1}.

Now we turn to the exponential decay estimate for the derivatives of $v_{1}$ and $v_{2}$, let $\partial$ denote any one of the two partial derivatives $\partial_{1}$ and $\partial_{2}$. Then we differentiate \eqref{4.38} to get
\be\label{4.45}
\triangle\left(\partial v_{j}\right)=\sum_{k=1}^{2} m_{jk}\mathrm{e}^{2\sum\limits_{l=1}^{2} b_{kl}^{-1}v_{k}}\left(2\sum_{l=1}^{2} b_{kl}^{-1}\right)\left(\partial v_{k}\right),~~~j,l=1,2.
\ee
Introduce the notation
\be
I=\mathrm{diag}\{1,1\},~~E(x)=\mathrm{diag}\left\{\mathrm{e}^{\frac{2\beta}{2\alpha-1}v_{1}(x)},\mathrm{e}^{v_{2}(x)}\right\},~~\mathbf{h}=\left(\partial v_{1},\partial v_{2}\right)^{\tau},\nn\\
\ee
then we see that \eqref{4.45} takes the matrix form
\ber\label{4.46}
&&\triangle\mathbf{h}=2MB^{-1}\mathbf{h}+2MB^{-1}\left(E(x)-I\right)\mathbf{h}\nn\\
&&~~~~\,\triangleq2D\mathbf{h}+2D\left(E(x)-I\right)\mathbf{h}.
\eer
Noting that there exists an invertible matrix
\be
T=\left(
\begin{matrix}
\dfrac{2\alpha-1}{2\left(\lambda_{3}-\alpha\right)}&~~~1\\[8pt]
1&~~~
\dfrac{2\left(\lambda_{4}-\alpha\right)}{2\alpha-1}
\end{matrix}
\right)\nn
\ee
such that $TDT^{-1}=\Lambda=\mathrm{diag}\{\lambda_{3},\lambda_{4}\}$, where $\alpha=\dfrac{3}{2}-\dfrac{1}{2N}$ and
\be\label{4.46b}
\lambda_{3}=\frac{2N+1+2\sqrt{N^{2}-N+\frac{1}{4}}}{4},~~~\lambda_{4}=\frac{2N+1-2\sqrt{N^{2}-N+\frac{1}{4}}}{4}
\ee
are two positive eigenvalues of the matrix $T$. Setting $\mathbf{H}=T\mathbf{h}$, then \eqref{4.46} can be rewritten as
\ber\label{4.46a}
&&\triangle\mathbf{H}=2\Lambda\mathbf{H}+2TD\left(E(x)-I\right)T^{-1}\mathbf{H}.
\eer
By straightforward calculations we have
\ber\label{4.47}
&&\triangle\mathbf{|H|}^{2}\geq2\mathbf{H}^{\tau}\triangle\mathbf{H}\nn\\
&&~~~~~~=4\mathbf{H}^{\tau}\Lambda\mathbf{H}+4\mathbf{H}^{\tau}TD\left(E(x)-I\right)T^{-1}\mathbf{H}\nn\\
&&~~~~~~\geq\lambda\mathbf{|H|}^{2}-b(x)\mathbf{|H|}^{2},~~x\in\mathbb{R}^{2}\backslash(0,0),
\eer
where
\be\label{4.48}
\lambda=\min\left\{\lambda_{3},\lambda_{4}\right\}>0
\ee
and the function $b(x)$ satisfies $b(x)\rightarrow0$ as $|x|\rightarrow\infty$. Consequently, as before, we infer that there exists a constant $C(\varepsilon)>0$ such that
\be\label{4.49}
\mathbf{|H|}^{2}\leq C(\varepsilon)\mathrm{e}^{-(1-\varepsilon)\sqrt{\lambda}|x|},~~\text{as}~~|x|\to \infty,
\ee
which gets estimate \eqref{3.19} for the derivatives of solutions.

We next are in a position to calculate the quantized integrals applying the decay estimates. Indeed, we have obtained the unique solution $(u_{1},u_{2})$ of the system \eqref{3.15} subject to the boundary condition \eqref{3.13} which vanish at infinity exponentially fast. Additionally, we deduce that $|\nabla(mu_{1}+2u_{2})|, |\nabla(pu_{1}+qu_{2})|$ all vanish at infinity at least as fast as $|x|^{-3}$, where $m, p, q$ are as defined by \eqref{3.19a}. Furthermore, from \eqref{4.0a}--\eqref{4.0b} and the exponential decay property of $|\nabla(mu_{1}+2u_{2})|$ and $|\nabla(pu_{1}+qu_{2})|$, we can see that $|\nabla(mP_{1}+2P_{2})|=O(|x|^{-3})$ and $|\nabla(pP_{1}+qP_{2})|=O(|x|^{-3})$ at infinity. Thus, according to the divergence theorem, we see that
\be\label{4.50}
\int_{\mathbb{R}^{2}} \triangle \left(mP_{1}(x)+2P_{2}(x)\right)\dd x=\int_{\mathbb{R}^{2}} \triangle \left(pP_{1}(x)+qP_{2}(x)\right)\dd x=0.
\ee
Now according to \eqref{4.0c}, \eqref{4.0d}, and \eqref{4.50} and the definitions of $\varphi_{i}(i=1,2)$, we can get the quantized integrals \eqref{3.20} as stated in Theorem \ref{th3.1}.


\begin{thebibliography}{99}

\bibitem{Abel}
T. Abel, A. Stebbins, P. Anninos and M.L. Norman, First structure formation. II. Cosmic string plus hot dark matter models, \emph{Astrophy. J.} {\bf 508} (1998) 530--534.

\bibitem{Abrokosov}
A.A. Abrikosov, On the magnetic properties of superconductors of the second group, \emph{Sov. Phys. JETP} {\bf 5} (1957) 1174--1182.

\bibitem{Actor}
A. Actor, Classical solutions of $SU(2)$ Yang--Mills theories, \emph{Rev. Mod. Phys.} {\bf 51} (1979) 461--525.

\bibitem{Ambjorn}
J. Ambj{\o}rn and P. Olesen, A condensate solution of the electroweak theory which interpolates between the broken and the symmetric phase, \emph{Nucl. Phys. B} {\bf 330} (1990) 193--204.

\bibitem{Bartolucci}
D. Bartolucci and G. Tarantello, Liouville type equations with singular data and their applications to periodic multivortices for the electroweak theory, \emph{Commun. Math. Phys.} {\bf 229} (2002) 3--47.

\bibitem{Bezryadina}
A. Bezryadina, E. Eugenieva and Z. Chen, Self--trapping and flipping of double--charged vortices in optically induced photonic lattices, \emph{Optics Lett.} {\bf 31} (2006) 2456--2458.

\bibitem{Bo}
S. Bolognesi, C. Chatterjee, J. Evslin, K. Konishi, K. Ohashi and L. Seveso, Geometry and dynamics of a coupled 4D--2D quantum field theory, \emph{J. High Energy Phys.} {\bf 01} (2016) 075.

\bibitem{Chen0}
S. Chen, X. Han, G. Lozano and F.A. Schaposnik, Existence theorems for non--Abelian Chern--Simons--Higgs vortices with flavor, \emph{J. Differential Equations} {\bf 259} (2015) 2458--2498.

\bibitem{Chen1}
S. Chen and Y. Yang, Existence of multiple vortices in supersymmetric gauge field theory, \emph{Proc. R. Soc. A} {\bf 468} (2012) 3923--3946.

\bibitem{Chen2}
S. Chen, R. Zhang and M. Zhu, Multiple vortices in the Aharony--Bergman--Jafferis--Maldacena model, \emph{Ann. Henri Poincar$\acute{e}$} {\bf 14} (2013) 1169--1192.

\bibitem{Fujimoto}
Y. Fujimoto and M. Nitta, Topological confinement of vortices in two--flavor dense QCD, \emph{J. High Energy Phys.} {\bf 09} (2021) 192.

\bibitem{Gilbarg}
D. Gilbarg and N. Trudinger, \emph{Elliptic Partial Differential Equations of Second Order}, Springer, Berlin and New York, 1977.

\bibitem{Hanany}
A. Hanany and D. Tong, Vortices, instantons and branes, \emph{J. High Energy Phys.} {\bf 07} (2003) 037.

\bibitem{Hartnoll}
S.A. Hartnoll, C.P. Herzog and G.T. Horowitz, Holographic superconductors, \emph{J. High Energy Phys.} {\bf 12} (2008) 015.

\bibitem{Hindmarsh}
M.B. Hindmarsh and T.W.B. Kibble, Cosmic strings, \emph{Rep. Prog. Phys.} {\bf 58} (1995) 477--562.

\bibitem{Jaffe}
A. Jaffe and C.H. Taubes, \emph{Vortices and Monopoles}, Birkh$\ddot{a}$user, Boston, 1980.

\bibitem{Julia}
B. Julia and A. Zee, Poles with both magnetic and electric charges in non--Abelian gauge theory, \emph{Phys. Rev. D} {\bf 11} (1975) 2227--2232.

\bibitem{Kawaguchi}
Y. Kawaguchi and T. Ohmi, Splitting instability of a multiply charged vortex in a Bose--Einstein condensate, \emph{Phys. Rev. A} {\bf 70} (2004) 043610.

\bibitem{Khomskii}
D.I. Khomskii and A. Freimuth, Charged vortices in high temperature superconductors, \emph{Phys. Rev. Lett.} {\bf 75} (1995) 1384--1386.

\bibitem{Lieb}
E.H. Lieb and Y. Yang, Non--Abelian vortices in supersymmetric gauge field theory via direct methods, \emph{Commun. Math. Phys.} {\bf 313} (2012) 445--478.

\bibitem{Mandelstam}
S. Mandelstam, Vortices and quark confinement in non--Abelian gauge theories, \emph{Phys. Lett. B} {\bf 53} (1975) 476--478.

\bibitem{Manton}
N. Manton and P. Sutcliffe, \emph{Topological Solitons}, Cambridge University Press, Cambridge, UK, 2004.

\bibitem{Nielsen}
H.B. Nielsen and P. Olesen, Vortex--line models for dual strings, \emph{Nucl. Phys. B} {\bf 61} (1973) 45--61.

\bibitem{Sakai}
N. Sakai and Y. Yang, Moduli space of BPS walls in supersymmetric gauge theories, \emph{Comm. Math. Phys.} {\bf 267} (2006) 783--800.

\bibitem{Shevchenko}
S.I. Shevchenko, Charged vortices in superfluid systems with pairing of spatially separated carriers, \emph{Phys. Rev. B} {\bf 67} (2003) 214515.

\bibitem{Shifman}
M. Shifman and A. Yung, \emph{Supersymmetric Solitons}, Cambridge University Press, Cambridge, UK, 2009.

\bibitem{Sokoloff}
J.B. Sokoloff, Charged vortex excitations in quantum Hall systems, \emph{Phys. Rev. B} {\bf 31} (1985) 1924--1928.

\bibitem{Stoer}
J. Stoer and R. Bulirsch, \emph{Introduction to Numerical Analysis}, Springer, New York, 1983.

\bibitem{Tallarita}
G. Tallarita, Non--Abelian vortices in holographic superconductors,\emph{Phys. Rev. D} {\bf93} (2016) 066011.

\bibitem{Tong}
D. Tong, Quantum vortex strings: a review, \emph{Ann. Phys.} {\bf 324} (2009) 30--52.

\bibitem{Vilenkin}
A. Vilenkin and E.P.S. Shellard, \emph{Cosmic Strings and Other Topological Defects}, Cambridge University Press, Cambridge, UK, 1994.

\bibitem{Yang}
Y. Yang, On a system of nonlinear elliptic equations arising in theoretical physics, \emph{J. Funct. Anal} {\bf 170} (2000) 1--36.

\bibitem{YangBook}
Y. Yang, \emph{Solitons in Field Theory and Nonlinear Analysis}, Springer, New York, 2001.

\bibitem{YM}
C.N. Yang and R.L. Mills, Conservation of isotopic spin and isotopic gauge invariance, \emph{Phys. Rev.} {\bf96} (1954) 191--195.

\end{thebibliography}
\end{document}